\documentclass[a4paper, 12pt]{article}
\usepackage{amssymb}
\usepackage{amsmath,amsthm}
\usepackage{graphicx}

\begin{document}

\title{Information propagation in a novel hierarchical network }
\author{Feng Fu, Lianghuan Liu, Long Wang\\
Laboratory for Intelligent Control,
Center for Systems and Control,\\
Department of Mechanics and Engineering Science,\\
Peking University, Beijing 100871, P. R. China}
\date{}
\maketitle

\begin{abstract}
A novel hierarchical model of complex network is proposed, and
information propagation process taking place on top of the network
is investigated. Our model is demonstrated to have small world
property. We found that the frequency distribution of refractory
elements number is bimodal and the location of initially chosen
seed infection affects the spreading range of the information.\\


\end{abstract}

PACS number(s): 89.75.Hc, 87.23.Kg

\section{Introduction}
The research in the realm of complex networks has been boosted
spectacularly these years. The beauty of complex networks
fascinates many scientists from different background to contribute
to this issue
\cite{strotagz_nature,Newman,Albert_Barabasi_review_2002,durrett_graph,Boccaletti_2006}.
Since the existence of random graph theory has been the important
transition from regular lattices to irregular networks, the
small-world (Watts-Strogatz) model and scale-free
(Barb\'{a}si-Albert) model have been the catalyst for the rise of
interest of the whole scientific world. Meanwhile, the dynamical
evolution of the networks (i.e. epidemic spreading in the
networks, cascading failure of the power grid, cooperation in the
networks and etc) has been also investigated. To our best
knowledge, although the Watts-Strogatz (WS) and Barb\'{a}si-Albert
(BA) model could reflect some properties of real complex networks,
both of them are not sufficient to elucidate the real networks.
Recently, the generic features of hierarchical structure and
self-similarity of complex networks shed light on development of
novel models \cite{social_hierarchy_science}.

The literature related to networks is like booming mushroom after
rain these years. Based upon Erd\H{o}s and R\'{e}nyi's prominent
work of random graph theory \cite{ER_1,ER_2}, Watts and Strogatz
first published their results about small-world phenomenon
\cite{Watts-Strogatz}, which is first discovered by the
sociologist Milgram \cite{milgram}. And almost at the next year,
Barab\'{a}si and Albert unravelled the twist of real complex
networks, such as World Wide Web, internet, scientific
collaboration networks and etc. They discovered the scale-free
topology of the real networks, namely, the degree distribution has
power law tail---$p(k)\sim k^{-\gamma}$. They also proposed an
algorithm to generate such scale-free networks which composed of
two factors: growth and preferential attachment
\cite{BA_science_1999}. After Barab\'{a}si's work, Krapivsky
studied a class of growing random networks and found that a power
law distribution of degree arises only for linear connection
kernel \cite{KSF_random}. Dynamical evolving processes in networks
such as age, fitness, epidemic spreading, rumor, cooperation and
so forth have also been investigated and some important results
are obtained
\cite{Watts-book-small-world,Kleinberg,P-S-V-critical-1,P-S-V-critical-2}.

\begin{figure}
\centering
\includegraphics[scale=0.7]{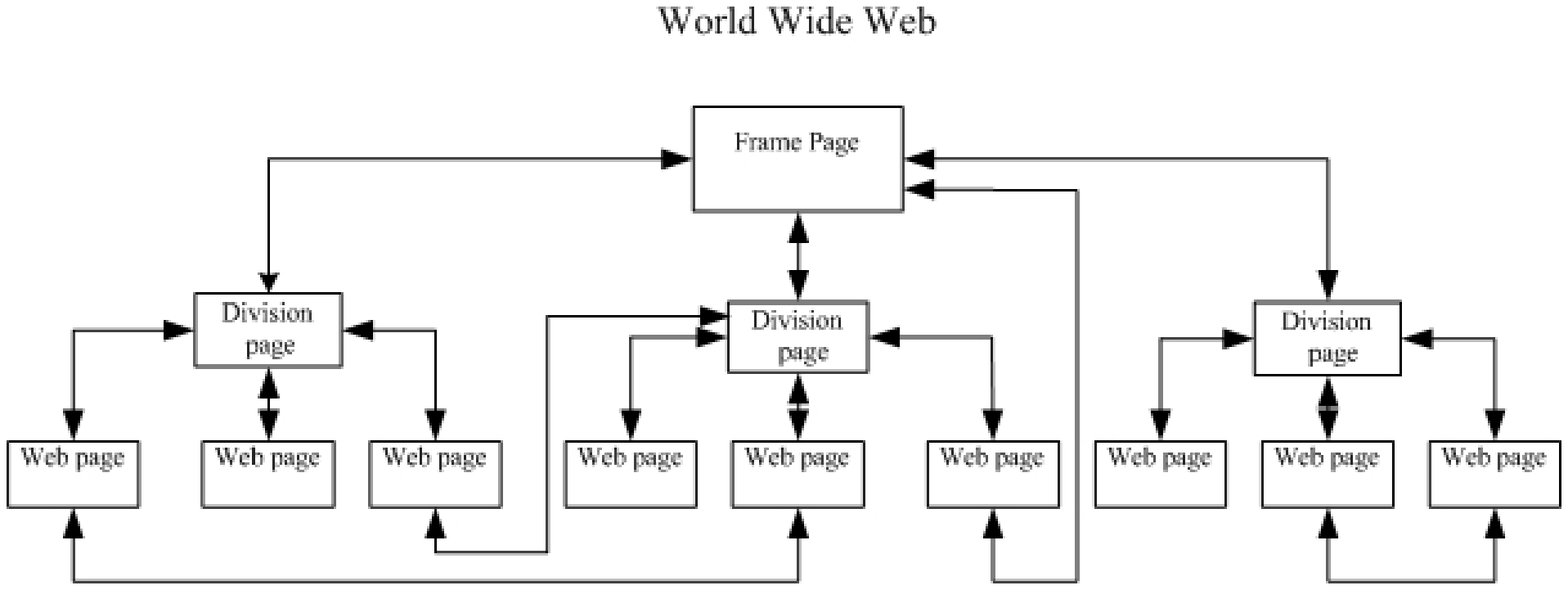}
\renewcommand{\figurename}{Fig.}
\caption{Hierarchical organizations in World Wide Web, which
denote communities of web pages with common interests.}
\label{webpage}
\end{figure}

Fundamental characteristics of real complex networks---small
worlds and scale-invariant topology, are essential to existent
network models. However, real networks are deviated from smooth
simulation results obtained from these ideal models. For instance,
the degree distribution of BA model obeys power law strictly,
whereas most of real networks only have power law tails.
Therefore, such discrepancy between models and empirical
measurements is attributable to previously disregarded, yet
generic feature of many real networks: their hierarchical
topology. Along with hierarchy, most networks are modular: many
small modules of highly connected nodes combine in a hierarchical
manner into larger, less cohesive units, and modules build up into
larger blocks and so on. In this way, modules are self-organized
into hierarchical structures. Indeed, in the World Wide Web, web
pages most link to the ones shared their common interests, and
these interdependent communities develop into larger groups
focusing on more general topics. In this sense, the WWW has its
``core'' which is composed of the most influential web sites, and
this core generates plenty of subdivisions and so on(see Fig.
\ref{webpage}). Accordingly, the WWW develops into the most huge
networks in the world with hierarchical structures.

Although hierarchical structures of complex networks are
investigated by many researchers, existent models which account
for hierarchy are very few \cite{Kima_2005}. Barab\'{a}si and his
collaborators constructed a model of network which has
hierarchical structure by iteration
\cite{hierarchy_arxiv,d-s-f_phsA}. Cellular metabolic networks are
also demonstrated to have hierarchical organizations
\cite{Gerdes,influx-meta}. However, hierarchical organizations are
ubiquitous in real complex networks, for instance, social
organizations or departments, the Internet at the domain level,
hierarchical architecture of web pages, and etc. Recently, real
networks are revealed that they commonly have fractal geometry
patterns corresponding to tree and snowflakes, namely,
self-similarity \cite{self-similarity_nature}. Consequently, it is
necessary and meaningful to depict the complex networks from their
hierarchical organization. Thus, this paper mainly deals with
hierarchy, and proposes an alternative generating approach of
complex networks.

Research on epidemic in networks has received increasing attention
since the inspiring work of Sudbury in random graph
\cite{Sudbury_1985}. The problem of epidemic deals with whether an
initially localized seed infection can spread to a substantial
part of the network. The discoveries of small-world and scale-free
networks model boost the studies of epidemic in these complex
networks. The pioneering work of Pastor-Satorras and Vespignani
found that there does not exist intrinsic spreading threshold for
scale-free networks \cite{Vespiganani_sf_epidemic}. By contrast,
there exits nonzero spreading threshold for the homogeneous and
small world networks \cite{Marro_book,Murray_biomath}. The
susceptible-infected-refractory(SIR) model has been widely used to
understand the spreading process of virus(for example,
transmission of influenza), rumor spreading and so forth
\cite{Masuda_2004,Zanette_smallworld_2001,Zanette_2002,Moreno_2004,Dodds_Watts_PRL_2004}.
Nevertheless, sometimes it is desirable to spread the ``epidemic''
quickly and effectively. Some typical examples are rumor-based
protocols for data dissemination and resource discovery on the
Internet and marketing campaigns using rumor-like strategies
(viral marketing) \cite{Moreno_2004_R}. Epidemic spreading process
is also investigated under the different hierarchical network
topologies \cite{Grabowski_2004,Zheng_2005,Huang_2006}.

In this paper, the information (rumor, news and facts)propagation
on top of our hierarchical network is investigated by using SIR
model. Once the ignorant gets a piece of information, there is a
high possibility that he/she would like to spread the information
to his/her neighbors in the network. By convenience, the spreading
rate is set to be unit, which is the possibility that a
susceptible individual is infected when contacted.

Our paper is organized as follows. Section II proposes the
hierarchical network model and construction procedure; Section III
discusses the SIR dynamics of information propagation on top of
the network; Section IV gives out the simulation results, and some
important conclusions are made in Section V.
\section{Hierarchical Network Model and Construction Procedure}

\begin{figure}
\centering
\includegraphics[scale=0.5]{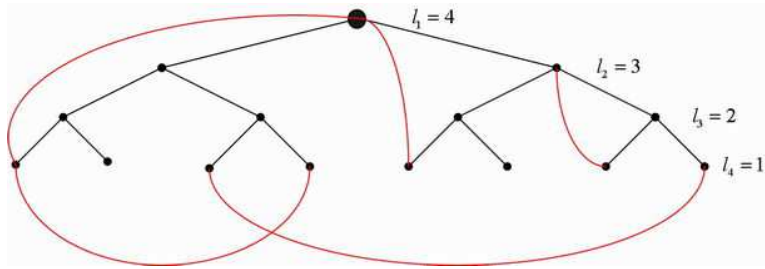}
\renewcommand{\figurename}{Fig.}
\caption{Illustration of typical hierarchical structures in
complex networks, and the red lines are the long-range connections
between nodes at different/same layers.} \label{hierarchy}
\end{figure}

Considering the undirected graph, we start with $m_0$ nodes at the
first layer, and each node connects to his distinct sons (not more
than $m$) which belong to the next layer and repeat this procedure
$l$ times, we obtain the skeleton. Such tree-like web is
hierarchical, hub-cascading and not sufficient to reflect the
complexity of real networks. However, despite the local links,
additional long-range direct connections (shortcuts) between nodes
should be necessary in order to make the average path short(see
Fig.\ref{hierarchy}). So for nodes belonging to the same layer,
connect to each other with possibility $p_{ij} = c_1e^{-\alpha
x_{ij}}$, where $\alpha > 0$ is tunable parameter and a measure of
homophily---the tendency of like to associate with like and $c_1 >
0$ is a normalizing constant. $x_{ij}$ is ``social distance'' of
node $i, j$, i.e. the lowest height of their common ancestor. When
$e^{-\alpha}\ll 1$, individuals will connect only to those most
similar to themselves. Besides, nodes of lower layer would prefer
to link to nodes of upper layers. Define the height of layers from
top to bottom as $l_j,j=1,2,\ldots,l$ respectively.  Node $k$ of
$i$th layer links to node $n$ of $j$th layer with possibility
$p_{kn} = \frac{l_j}{\sum_k l_k}e^{-\alpha x_{kn}}$, ($l_j
> l_i$). Given such connection kernel, individuals prefer to link to the nodes of
higher layers and most similar to themselves. After addition of
these shortcuts, a hierarchical network topology is obtained.

\section{Dynamics of Information Propagation}
Considering the information propagation process taking place on
our hierarchical social networks, assume each node has three
different states, analogously to SIR model in epidemiology:
susceptible (S), infected (I), refractory (R). Node $i$ could only
be in one of the three possible states at each time step of
evolution. At the beginning, only one individual randomly chosen
as seed obtains a piece of information and becomes infected. And
all the remnant population is susceptible. The evolution rules act
as follows \cite{Frauenthal_1980}. At each time step, an
individual $i$ is chosen at random from the infected population.
This individual contacts one of his/her neighbors $j$. If $j$ is
in the susceptible state, and it becomes infected. Otherwise, if
$j$ is already infected or refractory, then $i$ becomes
refractory. This progress continues until there are no longer
infected individuals in the population. The time the whole process
takes is the information lifetime $T$. For $t \geq T$, the number
of refractory nodes $N_R$ represents the number of nodes that are
infected in the network and acquire the information finally.

Qualitatively speaking, the dynamical process of information
propagation could be summarized as follows. In the early stage of
the evolution, the number of infected nodes increases and with
lower rate, the refractory population grows as well. Consequently,
the increasing contacts of infected nodes between themselves and
with refractory ones results in the decline of infected
population. At the end, the infected individuals vanish, and the
population consists of a group of $N_R$ refractory individuals,
who have been infected during the evolution and a group of
susceptible ones who have never heard the information. In
practice, the distribution $f(N_R)$ of the number $N_R$ of
refractory nodes is usually investigated over a large number of
realizations of the system evolution. Thus the average fraction
$r$ of the refractory population could be calculated as $r=\langle
N_R \rangle/N=N^{-1}\sum_{N_R=0}^NN_Rf(N_R)$. This could be
measurement of the influence of the information propagation on top
of complex networks.

For homogenous networks (e.g. random networks), the analytical
mathematical theory of SIR dynamics with unit spread rate suggests
that the fraction $r$ approaches to a limit $r^*$ as the system
size $N \to \infty$, where $r^*$ is the nontrivial solution of the
transcendental equation: $r^*=1-e^{-2r^*}$, i.e.
$r^*\approx0.796$\cite{Sudbury_1985}. In other words, about $20\%$
of the population would never have the chances to receive the
information. For more complex networks, we have to rely on the
numerical simulations to investigate the temporal behavior of the
information propagation process.

\section{Simulation Results}
Statistical quantities of such hierarchical network are computed,
namely, degree distribution, clustering coefficient and average
path length.

In our simulations, $c_1 = 0.2, \alpha = 0.5, m = 6$, and the
number of layers depends upon generating process, different
numbers of nodes $N =1000,2000,$
$3000,4000,5000,6000,7000,8000,9000,10000,15000,20000
,25000,30000$ are computed respectively.

Our algorithm is as follows:

\begin{itemize}
\item{step 1: $n_1$ nodes are selected as members of the first
layer, where $n_1$ is a random number not more than
$\frac{1}{100}N$ ;} \item{step 2: generate a random number $r$,
which is the number of nodes of next layer. If $r > n_i^m$ or $r >
N_{left}$, let $r=\min\{n_i^m, N_{left}\}$, where $n_i$ is the
number of node of upper layer $i$, and $N_{left}$ is the number of
unallocated nodes. And then stochastically find these $r$ nodes'
father in upper layer until every node finds their father, and
every father has no more than $m$ sons;}\item{step 3: repeat this
procedure until all nodes are allocated to corresponding
layers;}\item{step 4: for nodes belonging to the same layer, each
pair of individuals connect to each other with possibility $p_{ij}
= c_1e^{-\alpha x_{ij}}$;}\item{step 5: begin with nodes of the
lowest layer, generate connections to nodes of upper layers with
possibility $p_{kn} = \frac{l_j}{\sum_k l_k}e^{-\alpha x_{kn}}$
until nodes of every two layers have been linked with the
possibility function.}
\end{itemize}

Clustering coefficient $C(N)$ is defined as the average clustering
coefficient $C_i$ of each node $i$, and $C_i =
\frac{2E_i}{k_i(k_i-1)}$, where $k_i$ is the degree of node $i$,
and $E_i$ is the number of existent links between its $k_i$
nearest neighbors. Average path length is the average length of
all existent connections between pairs of nodes. Diameter of the
graph is the maximum path length of pairs of connected nodes. The
following table is the summary of simulation results.
\begin{tabular}{|c|c|c|c|}\hline
\multicolumn{4}{|c|}{\bfseries Statistical Quantities of Different
System Size}\\
\hline Number of Nodes & Clustering Coefficient & Average Path
Length &
Diameter\\[0.5ex]
\hline
1000 &   0.297 &   4.2 & 11 \\
\hline
2000 & 0.180 & 3.1 & 11\\
\hline
3000 & 0.191  &  4.3 & 12 \\
\hline
4000  &  0.153 & 3.5 &   9 \\
\hline
5000 & 0.131 &   3.6 & 12\\
\hline
6000 &   0.130 & 3.4& 10 \\
\hline
7000 & 0.137 &   3.4 & 18 \\
\hline
8000 & 0.116 & 2.9 & 12 \\
\hline
9000 &  0.145 & 2.7 &   15\\
\hline
10000 & 0.135 & 3.2 & 16\\
\hline
15000 & 0.100 & 3.5 & 15\\
\hline
20000 & 0.137 & 3.9 & 17\\
\hline
25000 & 0.079 & 3.0 &11\\
\hline
30000 & 0.083 & 2.7 & 11\\
\hline
\end{tabular}\\

Our simulation results demonstrate that our hierarchical network
has small world property, as the average path length is
considerably short, compared to the same size of random graph.
Furthermore, the clustering coefficient is also relatively high,
which means that such network has small worlds phenomenon.
Meanwhile, the long-range connections between nodes are essential
to make the network searchable, i.e. short average path length.
Equivalently, our hierarchy network has good searchability.

\begin{figure}
\centering
\includegraphics[scale=0.1]{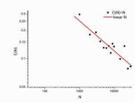}
\renewcommand{\figurename}{Fig.}
\caption{Plot of $C(N)$ vs. system size $N$ on logarithmic scale.}
\label{C_N}
\end{figure}

From the Fig.\ref{C_N}, our simulation results indicate that the
clustering coefficient depends on the system size as $C(N) \sim
N^{-0.32}$, which is significantly larger for large $N$ than the
random network prediction $C(N) \sim N^{-1}$. Indeed, many real
networks have the same property as our model. Yet, degree
distribution of our model does not follow the power law.

Concerning the information propagation process on top of our
hierarchical network, we make some slight variations of the
simulated network model in order to reveal more new phenomena
because of the hierarchical topology, other than singularity of
the structure. Firstly, the node number of the first layer is not
more than $1/1000$ of the total population. Secondly, exception
for the nodes of the lowest layer, every node has constant $m$
sons which belong to the next layer. And also the situation that
$m$ is subject to a uniform gaussian distribution is considered.
Thirdly, node $i$ connects to the other node $j$ belonging to the
same layer with probability $p_{ij}=c_1e^{-\alpha x_{ij}}$, and
links to his ancestor $j$ in upper layers with probability
$p_{ij}=\frac{l_{j}}{\sum_kl_k}$.

\begin{figure}
\centering
\includegraphics[scale=0.4]{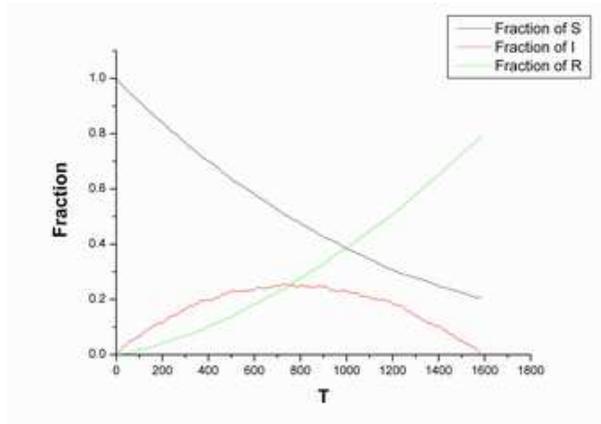}
\renewcommand{\figurename}{Fig.}
\caption{The dynamic evolution of information propagation. System
size $N=1000$.} \label{evolution_1k}
\end{figure}

The information propagation process was performed by using SIR
model as above instructions. The initial condition was that
randomly chosen individual got a piece of information and became
``infected'' (spreader). He/she would spread the information to
his/her ``susceptible'' (ignorant) neighbors until he/she became
``refractory'' (stifler). Fig.\ref{evolution_1k} shows the
information propagation process with population $N=1000$. The
fraction of infected individuals $i(t)$ increases faster than the
fraction of refractory individuals $r(t)$ in the beginning. Then
$i(t)$ decreases to zero at the end of evolution. The fraction of
refractory population is $79.4\%$ when the information propagation
process terminates. Accordingly, there are still $20.6\%$ of the
total population who never acquire the information through their
neighbors.

\begin{figure}
\centering
\includegraphics[scale=0.4]{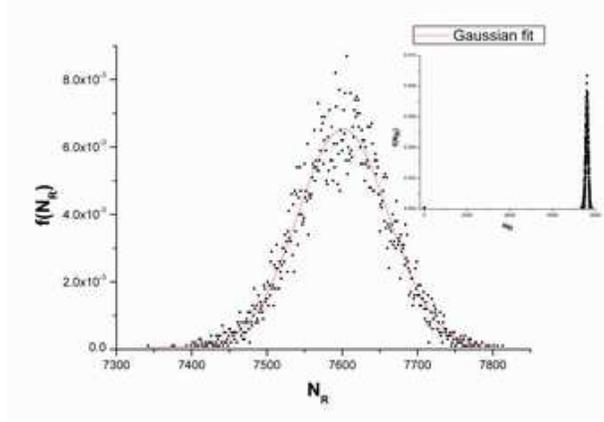}
\renewcommand{\figurename}{Fig.}
\caption{The insert shows frequency distribution $f(N_R)$ of $N_R$
over $10^4$ realizations. System size is $N=10000$, $m=32$.
Gaussian fit of the data shows that for large $N_R$, $f(N_R)$
approximately follows a gaussian distribution.}
\label{NR_distribution}
\end{figure}

\begin{figure}
\centering
\includegraphics[scale=0.1]{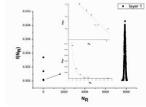}
\renewcommand{\figurename}{Fig.}
\caption{The frequency distribution $f(N_R)$ over $10^4$
realizations with initially infected individual belonging to the
first layer. System size $N=10000$,$m=8$. The lower insert shows
the $N_R$'s distribution near zero and the upper insert is the
corresponding plot on logarithmic scale , where the red line is
linear fit. } \label{N_R_layer1_constant}
\end{figure}

\begin{figure}
\centering
\includegraphics[scale=0.1]{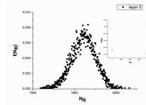}
\renewcommand{\figurename}{Fig.}
\caption{The large-$N_R$ structure. The insert shows the frequency
distribution $f(N_R)$ over $10^4$ realizations with initially
infected individuals belonging to the third layer. System size
$N=10000$,$m=8$.} \label{N_R_layer3_constant}
\end{figure}

One of the important practical aspects of information propagation
is whether or not it could reach a high number of individuals.
Hence, we have studied the distribution of the number $N_R$ of R
elements at the end of evolution. Fig.\ref{NR_distribution} shows
the normalized frequency $f(N_R)$ obtained from series of $10^4$
realizations for large $N_R$ with system size $N=10000$, $m=32$.
Most of the $N_R$ falls into the large-$N_R$ structure, whereas
very few $N_R$ is in the neighborhood of zero. In
Fig.\ref{NR_distribution}, only one point is adjacent to zero.
Therefore, we could conjecture that this particular distribution
$f(N_R)$ is attributable to the network's hierarchical topology.
To investigate more details about the effect of the network
topology on information propagation, different initial conditions
are considered. Different groups of realizations corresponding to
the layers are conducted as follows: for each group of
realizations, the initially selected at random infected individual
always belongs to the same layer. Figure \ref{N_R_layer1_constant}
shows the distribution $f(N_R)$ over $10^4$ realizations with the
initially infected node in the first layer. Obviously, here the
distribution $f(N_R)$ is bimodal, with a maximum close to $N_R=0$
and an additional bump for larger $N_R$. For small $N_R$ near
zero, the frequency follows a power law, $f(N_R)\sim
N_R^{-\alpha}$ ($\alpha\approx 3$). Contributions to this zone of
the distribution result from realizations where propagation ceases
before a shortcut is reached. Howbeit, $f(N_R)$ approximately
obeys gaussian distribution for large $N_R$. In a typical
realization contributing to this area of the distribution, large
quantity of contacts occur through shortcuts and a finite portion
of the population becomes infected. Our simulation results show
that the positions of initially randomly selected individual
significantly affect the spreading range of information. When
chosen individual belongs the lower layer, the distribution
$f(N_R)$ has more points falling into the large-$N_R$ structure
(see Fig. \ref{N_R_layer3_constant}, for the initial individuals
lying in the third layer, almost all of the points consist of the
gaussian distribution). The average fraction of refractory
individuals $r$ over large number of realizations is computed.
With system size $N=10000$, $m=8$, the hierarchical network
realization has 5 layers. The fractions $r$ corresponding to each
layer are $77.21\%, 77.91\%, 78.34\%, 78.27\%, 78.28\% $
respectively. We found that the discrepancies of the spreading
range are remarkable when initially selected nodes belong to
different layers. Thereby, the hierarchical topology intrinsically
affects the shape of distribution function $f(N_R)$(see
Fig.\ref{N_R_layer1_constant},\ref{N_R_layer3_constant} for
comparison).

\begin{figure}
\centering
\includegraphics[scale=0.4]{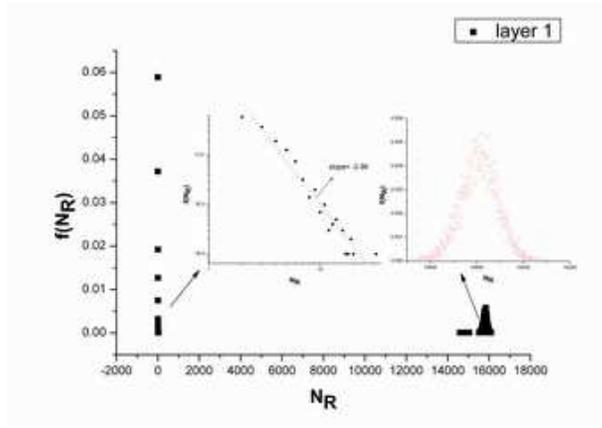}
\renewcommand{\figurename}{Fig.}
\caption{The plot of the situation that $m$ is subject to gaussian
distribution over $10^4$ realizations. The left insert shows the
detailed plot of the distribution near the zero on logarithmic
scale. The red line is linear fit of which slope is about -2.86.
And the right shows the large-$N_R$ structure which is
approximately a gaussian distribution. The system size is
$N=20000$. The mean value $\mu$ of $m$ is 8, and the standard
deviation $\delta$ is 2.  } \label{N_R_layer1_gaussian}
\end{figure}

Additionally, taking account for slight variation of the
hierarchical topology, i.e. $m$ is not constant for each node, but
subject to a uniform gaussian distribution. The reason why this
regime is adapted is to investigate whether or not irregularity of
the hierarchy would affect the information propagation.
Fig.\ref{N_R_layer1_gaussian} shows that the frequency
distribution of $N_R$ by using the initial condition that the
randomly chosen node to be infected lies in the first layer (In
this realization, there are total 8 layers). The distribution is
similar to the case where $m$ is constant. The distribution of
small $N_R$ near zero also follows power law (see
Fig.\ref{N_R_layer1_gaussian}, $\alpha \approx 2.86$), whereas the
large-$N_R$ structure approximately obeys gaussian distribution.
As well as the situation that $m$ is constant, along with the
initially infected individual falling into the lower layers, the
number of data points near the zero decreases. In other words, the
spreading range of the information increases as the initially
selected node lies in the lower layers. In this case, the average
fractions according to each layer are $67.52\%, 70.17\%, 74.99\%,
78.24\%, 78.94\%, 78.98\%, 78.98\%, 75.23\%$ respectively.

Our simulations results also demonstrate that the possible maximum
fraction of R individuals in the end of evolution depends upon the
position of the initially infected individual in the hierarchical
network. When the initially selected individual belongs to a
certain intermediate layer, the maximum spreading range of the
information would be achieved.

\begin{figure}
\centering
\includegraphics[scale=0.35]{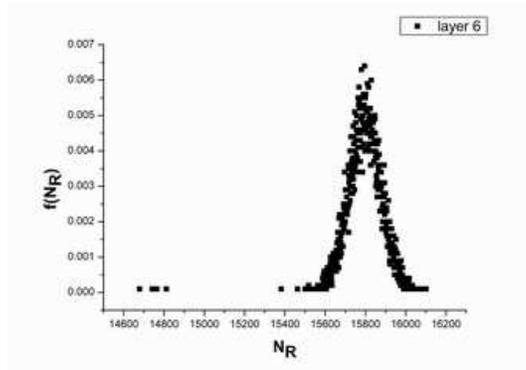}
\caption{Plot of the $N_R$ distribution with the initial infection
lying in the 6th layer. The other initial conditions are the same
as Fig. \ref{N_R_layer1_gaussian}.} \label{N_R_layer6_gaussian}
\end{figure}

\section{Conclusion}
In conclusion, the top layer of our model could be illustrated as
the core of real networks, for instance, the most popular web
sites. The model we propose in this paper has nontrivial
properties: it has short average path length and high clustering
coefficient, i.e. small world property; average clustering
coefficient depends on system size. From these two respects, our
preliminary model might be enlightening to reveal mechanism of
complex networks, and reflects characteristics of real networks.
Besides, the information propagation on top of our hierarchical
network is investigated by using SIR model. It is observed that
distribution $f(N_R)$ is bimodal: for small $N_R$ near the zero,
it follows power law whereas for large $N_R$, it is approximately
a gaussian distribution. We also found that the spreading range
approaches to maximum when the initially chosen individual is in a
certain intermediate layer. Our simulation results indicate that
the optimal strategy to spread the information efficiently is to
make the individuals belonging to the intermediate layers obtain
the information and would like to spread it to her/his neighbors.
It is also helpful to develop methods to control the information
propagation in networks. For example, to prevent the rumor from
mongering. Besides, our results demonstrate for a generic
realization, the fraction of population who obtain the information
is less than the limit of fraction of R individuals. Namely, there
are always ignorant individuals who are not aware of the
information. As a result, it is also significant to develop novel
propagandistic strategy to make almost everyone acquire the
information.


\end{document}